\documentclass{article}
\usepackage{graphicx}
\usepackage{amsmath}
\usepackage{amssymb}

\def\Z{\hbox{\bf Z}}

\def\R{\hbox{\bf R}}

\let\ds=\displaystyle

\makeatletter
\renewcommand{\section}{\@startsection
    {section}{1}{0mm}{.3ex}
    {.1ex}{\bf}}
\makeatother

\begin{document}

\title{\Large\bf SINC APPROXIMATION OF THE HEAT DISTRIBUTION ON THE BOUNDARY
 OF A TWO-DIMENSIONAL FINITE SLAB\thanks {Supported by the Council for Natural
 Sciences of Vietnam}}%
\author{\bf P. N. Dinh Alain\thanks{Mathematics Department, Mapmo UMR 6628,
 BP 67-59, Orleans Cedex 2, France. Email : alain.pham@math.cnrs.fr} - P. H. Quan\thanks{HoChiMinh
City National University, Department of Mathematics and Informatics, 227
Nguyen Van Cu, Q. 5, HoChiMinh City, VietNam. Email : tquan@pmail.vnn.vn}
\addtocounter{footnote}{-1}
- D. D. Trong\footnotemark}
\date{}
\maketitle%

{\bf Abstract :} We consider the two-dimensional problem of recovering
globally in time the heat distribution on the surface of a layer inside of a
heat conducting body from two interior temperature measurements. The problem
 is ill-posed. The approximation function is represented by a two-dimensional
 Sinc series and the error estimate is given.


{\bf Mathematics Subjects Classification 2000. 35K05, 31A25, 44A35}
\\
\\

{\bf 1. Introduction}

In this paper, we consider the problem of recovering the heat distribution on
 the surface of a thin layer inside of a heat conducting body from transient
temperature measurements. The problem is raised in many applications in
Physics and Geology. In fact, in many physical situation (see, e.g. [B]) we
 cannot attach a temperature sensor at the surface of the body (for example,
the skin of a missile). On the other hand, we can easily measure the
temperature history at an interior point of the body. Hence, to get the
heating history in the body, we want to use temperature measured in the
interior of the body. In the simplest model, the heat-conducting body is
assumed to have a constant conductivity and represented by the half-line $x >
 0$ (see, e.g. [C, EM, LN, TV]),. While giving many useful results, this
model is not suitable for the case of a body having a series of superposed
layers, each of which has a constant conductivity.

Precisely, we shall consider the problem corresponding to a thin layer of the
 body represented by the strip $\R \times (0,2)$, say. Let $u$ be the
 temperature in the strip. For the uniqueness of solution, we shall have to
measure the temperature history at two interior lines $\R \times \{y =
 1\}$ and $\R \times \{y = 2\}$. From these measurements, we can
identify uniquely the heating history inside of the layer (see, e.g., [B]).
 The problem is of finding the surface heat distribution $u(x,0,t) = v(x,t)$.
 In fact, despite uniqueness, the global solution in $L^2(\R \times
\R_+)$ is unstability and hence, in this point of view, a sort of
 regularization is in order.

As discussed in the latter paragraph, the main purpose of our paper is to
present a regularization of the problem. Moreover, an effective way of
approximating the heat distribution $v$ is also worthy of considering.
There are many methods for regularizing the equation (see [TA, B]). In the
most common scheme (see, e.g., [Blackwell]), the computation is divided into
 two steps. In the first step, one considers the problem of finding the heat
 flux $u_y(x,1,t)$ from the interior measurements $u(x,1,t), u(x,2,t)$. The
problem is classical and can be changed to the one of finding the solution of
 a convolution of Volterra type which can be solved in any finite time
interval by the iteration (see, e.g., [F]). But even in the "classical"
problem as mentioned in many documents, it is worthy of insisting that the
problem is {\sl ill-posed} if we consider the problem over the whole time
interval $\R_+$ with respect to the $L^2$-norm and the literature on
 this way is very scarce. In the second step, one considers the "really" ill-
posed problem which is of recovering the surface temperature history
$u(x,0,t)$ from data $u(x,1,t), u_y(x,1,t)$. In the present paper, we can
(and shall) regularize the function $u(x,0,t)$ {\sl directly} from $u(x,1,t),
 u(x,2,t)$ {\sl without using} the flux function $u_y(x,1,t)$. We emphasize
that, using our method, we can unify two steps of the classical scheme and
find simultaneously two functions $u(x,0,t)$ and $u_y(x,1,t)$. However, our
main purpose is of regularizing the surface temperature. Hence, we shall omit
 the problem of finding the interior heat flux $u_y(x,1,t)$.

Moreover, using the method of truncated integration, one can approximate the
 Fourier transform of the solution by a function having a {\sl compact
support} in $\R^2$. Therefore, the solution can be represented by an
expansion of two dimensional Sinc series (see [AGTV]). The Sinc method is
based on the Cardinal functions
\begin{eqnarray}
S(p,d)(z)=\frac{sin[(\pi (z-pd))/d]}{\pi (z - pd)/d},\,\,\,p \in \Z,
 d > 0
\nonumber
\end{eqnarray}
which dates back to the works of many mathematicians (Bohr, de la Vallee
Poussin, E. T. Whittaker, ...). The one dimensional version of the method is
studied very clearly and completely in [S]. Some primary results related to
the two dimensional Sinc approximation were given in [AGTV]. As is known, the
 Sinc series converges very rapidly at an incredible $0(e^{-cn^{1/2}})$ rate,
 where $c > 0$ and $n$ is the dimension of approximation (see [S]). Hence,
this method, which is new in our knowledge, is very effective.

The remainder of the present paper is divided into two sections. In Section 2
, we state precisely the problem, change it into an integral equation of
convolution type, and state main results of our paper. In Section 3, we give
 the proof of main results.\\

{\bf 2. Notations and main results}

Consider the problem of determining the heat distribution
\begin{eqnarray}
u(x,0,t) = v(x,t)
\end{eqnarray}
where $u$ satisfies
\begin{eqnarray}
\Delta u - \frac{{\partial u}}{{\partial t}} = 0\,\,\,\,\, x \in \R,\,
\,0 < y < 2,\,\,t > 0,
\end{eqnarray}
subject to the boundary conditions
\begin{eqnarray}
u(x,2,t) = g(x,t),\,\,\,x \in \R,\,\,\,t > 0,
\end{eqnarray}
\begin{eqnarray}
u(x,1,t) = f(x,t),\,\,\,x \in \R,\,\,\,t > 0,
\end{eqnarray}
and the initial condition
\begin{eqnarray}
u(x,y,0) = 0,\,\,\,x \in \R,\,\,\,0 < y < 2.
\end{eqnarray}

Here $f, g$ are given. We shall transform the problem (1)-(5) into a
convolution intergral equation.

Put
\begin{eqnarray}
\Gamma (x,y,t,\xi ,\eta ,\tau ) = \frac{1}{{4\pi (t - \tau )}}\exp \left(
 { - \frac{{(x - \xi )^2  + (y - \eta )^2 }}{{4(t - \tau )}}} \right)
\end{eqnarray}

and
\begin{eqnarray}
G(x,y,t,\xi ,\eta ,\tau ) = \Gamma (x,y,t,\xi ,\eta ,\tau ) - \Gamma (x,4 -y,
t,\xi ,\eta ,\tau ).
\end{eqnarray}

We have
\begin{eqnarray}
G_{\xi \xi }  + G_{\eta \eta }  + G_\tau   = 0.
\nonumber
\end{eqnarray}

Integrating the identity
\begin{eqnarray}
div(G\nabla u - u\nabla G) - \frac{\partial }{{\partial \tau }}(uG) = 0
\nonumber
\end{eqnarray}
over the domain $\R \times (1,2) \times (0,t - \varepsilon)$ and
letting $\varepsilon  \to 0$, we have
\begin{eqnarray}
\int\limits_{ - \infty }^{ + \infty } {\int\limits_0^t {g(\xi ,\tau )G_\eta
 (x,y,t,\xi ,2,\tau )d\xi d\tau } } &+& \int\limits_{ - \infty }^{ + \infty }
 {\int\limits_0^t {G(x,y,t,\xi ,1,\tau )u_y(\xi, 1, \tau )d\xi d\tau } }
\nonumber\\
&-& \int\limits_{ - \infty }^{ + \infty } {\int\limits_0^t {f(\xi ,\tau
)G_\eta  (x,y,t,\xi ,1,\tau )d\xi d\tau } }  + u(x,y,t) = 0.
\nonumber
\end{eqnarray}

Hence
\begin{eqnarray}
&&\int\limits_{ - \infty }^{ + \infty } {\int\limits_0^t {G(x,y,t,\xi ,1,\tau
 )u_y(\xi, 1, \tau )d\xi d\tau } }  =  - u(x,y,t) +
\nonumber\\
 &&\int\limits_{ - \infty }^{ + \infty } {\int\limits_0^t {G_\eta  (x,y,t,\xi
 ,1,\tau )f(\xi ,\tau )d\xi d\tau } }  - \int\limits_{ - \infty }^{ + \infty}
 {\int\limits_0^t {g(\xi ,\tau )G_\eta  (x,y,t,\xi ,2,\tau )d\xi d\tau } }.
\end{eqnarray}

Letting $y \to 1^ +$ in (8), we have
\begin{eqnarray}
&&\int\limits_{ - \infty }^{ + \infty } {\int\limits_0^t {\left[ {\frac{1}
{{2\pi (t - \tau )}}\exp \left( { - \frac{{(x - \xi )^2 }}{{4(t - \tau )}}}
 \right) - \frac{1}{{2\pi (t - \tau )}}\exp \left( { - \frac{{(x - \xi )^2  +
 4}}{{4(t - \tau )}}} \right)} \right]u_y(\xi, 1, \tau )d\xi d\tau } }
\nonumber\\
&& =  - f(x,t) - \frac{1}{{2\pi }}\int\limits_{ - \infty }^{ + \infty }
 {\int\limits_0^t {\frac{1}{{(t - \tau )^2 }}\exp \left( { - \frac{{(x -
\xi)^2  + 4}}{{4(t - \tau )}}} \right)f(\xi ,\tau )d\xi d\tau } }
\nonumber\\
&&+ \frac{1}{{2\pi }}\int\limits_{ - \infty }^{ + \infty } {\int\limits_0^t
{g(\xi ,\tau )\frac{1}{{(t - \tau )^2 }}\exp \left( { - \frac{{(x - \xi )^2
 + 1}}{{4(t - \tau )}}} \right)d\xi d\tau } }.
\end{eqnarray}

We put $N(x,y,t,\xi ,\eta ,\tau ) = \Gamma (x,y,t,\xi ,\eta ,\tau) - \Gamma
(x,- y,t,\xi ,\eta ,\tau )$

Integrating the identity
 \begin{eqnarray}
div (N\nabla u - u \nabla N) - \frac{\partial}{\partial\tau}(uN) = 0
\nonumber
\end{eqnarray}
over the domain $(-n,n) \times (0,1) \times (0,t - \varepsilon)$ and letting
$n \to \infty, \varepsilon  \to 0$
\begin{eqnarray}
\int\limits_{ - \infty }^{ + \infty } {\int\limits_0^t {N(x,y,t,\xi,1,\tau
 )u_y(\xi,1,\tau)d\xi d\tau } } &-& \int\limits_{ - \infty }^{ + \infty }
{\int\limits_0^t {f(\xi,\tau)N_\eta(x,y,t,\xi ,1,\tau )d\xi d\tau } }
\nonumber\\
+ \int\limits_{ - \infty }^{ + \infty } {\int\limits_0^t {v(\xi ,\tau
 )N_\eta(x,y,t,\xi ,0,\tau)d\xi d\tau } }  &-& u(x,y,t) = 0.
\end{eqnarray}

Letting $y \to 1^ -$, the identity (10) becomes
\begin{eqnarray}
&&\frac{1}{{2\pi }}\int\limits_{ - \infty }^{ + \infty } {\int\limits_0^t
 {\frac{1}{{t - \tau }}\left[ {\exp \left( { - \frac{{(x - \xi )^2 }}{{4(t -
 \tau )}}} \right) - \exp \left( { - \frac{{(x - \xi )^2  + 4}}{{4(t - \tau
 )}}} \right)} \right]u_y(\xi,1,\tau )d\xi d\tau } }
\nonumber\\
&& - \frac{1}{{2\pi }}\int\limits_{ - \infty }^{ + \infty } {\int\limits_0^t
 {f(\xi ,\tau )\frac{1}{{(t - \tau )^2 }}\exp \left( { - \frac{{(x - \xi )^2
 + 4}}{{4(t - \tau )}}} \right)d\xi d\tau } }
\nonumber\\
&& + \frac{1}{{2\pi }}\int\limits_{ - \infty }^{ + \infty } {\int\limits_0^t
{\frac{1}{{(t - \tau )^2 }}\exp \left( { - \frac{{(x - \xi )^2  + 1}}{{4(t -
\tau )}}} \right)v(\xi ,\tau )d\xi d\tau } }  - 3f(x,t) = 0
\end{eqnarray}

From (9) and (11), we have the main convolution equation
\begin{eqnarray}
&&-\frac{1}{\pi }\int\limits_{ - \infty }^{ + \infty } {\int\limits_0^t
 {\frac{1}{{(t - \tau )^2 }}\exp \left( { - \frac{{(x - \xi )^2  + 4}}{{4(t -
 \tau )}}} \right)f(\xi ,\tau )d\xi d\tau } }
\nonumber\\
&& + \frac{1}{{2\pi }}\int\limits_{ - \infty }^{ + \infty } {\int\limits_0^t
 {g(\xi ,\tau )\frac{1}{{(t - \tau )^2 }}\exp \left( { - \frac{{(x - \xi )^2
 + 1}}{{4(t - \tau )}}} \right)d\xi d\tau } }
\nonumber\\
&& + \frac{1}{{2\pi }}\int\limits_{ - \infty }^{ + \infty } {\int\limits_0^t
 {\frac{1}{{(t - \tau )^2 }}\exp \left( { - \frac{{(x - \xi )^2  + 1}}{{4(t -
 \tau )}}} \right)v(\xi ,\tau )d\xi d\tau } }  - 4f(x,t) = 0
\nonumber
\end{eqnarray}
which can be rewritten as
\begin{eqnarray}
S * v(x,t) =  2R * f(x,t) - S * g(x,t) + 4f(x,t)
\end{eqnarray}
where we define that $v(x,t) = f(x,t) = g(x,t) = 0$ as $t<0$,
\begin{eqnarray}
R(x,t) = \left\{ \begin{array}{l}
 \frac{1}{{t^2 }}\exp \left( { - \frac{{x^2  + 4}}{{4t}}} \right)\,\,\,\,\,
\,\,\,\,\,\,(x,t) \in  \R \times [0, + \infty ) \\
 \,\,\,\,\,\,\,0\,\,\,\,\,\,\,\,\,\,\,\,\,\,\,\,\,\,\,\,\,\,\,\,\,\,\,\,\,\,
\,\,\,\,\,(x,t) \in  \R \times ( - \infty ,0) \\
 \end{array} \right.
\nonumber
\end{eqnarray}
and
\begin{eqnarray}
S(x,t) = \left\{ \begin{array}{l}
 \frac{1}{{t^2 }}\exp \left( { - \frac{{x^2  + 1}}{{4t}}} \right)\,\,\,\,\,
\,\,\,\,\,\,(x,t) \in \R \times [0, + \infty ) \\
 \,\,\,\,\,\,\,0\,\,\,\,\,\,\,\,\,\,\,\,\,\,\,\,\,\,\,\,\,\,\,\,\,\,\,\,\,\,
\,\,\,\,\,(x,t) \in \R \times ( - \infty ,0) \\
 \end{array} \right.
\nonumber
\end{eqnarray}

Put $F(x,t) =   2R * f(x,t) - S * g(x,t) + 4f(x,t)$.

Taking the Fourier-transform of (12), we get
\begin{eqnarray}
\hat S(z,r)\hat v(z,r) = \hat F(z,r)
\nonumber
\end{eqnarray}
where
\begin{eqnarray}
\hat S(z,r) &=& \frac{1}{{2\pi }}\int\limits_{ - \infty }^{ + \infty
}\int\limits_{ - \infty }^{ + \infty } {S(z,r)e^{ - i(xz + tr)} dxdt}
\nonumber\\
&=& 2e^{ - \frac{1}{{\sqrt 2 }}\sqrt {\sqrt {z^4  + r^2 }  + z^2 } } \left[
 {\cos \frac{1}{{\sqrt 2 }}\sqrt {\sqrt {z^4  + r^2 }  - z^2 }  - isgn(r)
\sin \frac{1}{{\sqrt 2 }}\sqrt {\sqrt {z^4  + r^2 }  - z^2 } } \right]
\nonumber
\end{eqnarray}
and
\begin{eqnarray}
\left| {\hat S(z,r)} \right| = 2e^{ - \frac{1}{{\sqrt 2 }}\sqrt {\sqrt {z^4
 + r^2 }  + z^2 } } .
\nonumber
\end{eqnarray}

We have

{\bf Theorem 1}
{\sl

Let $\gamma  \in (0,2)$ and $\varepsilon  \in \left( {0,e^{\frac{-3}{\gamma
}} } \right)$.

Assume that $v_0 \in L^2(\R^2)$ is the (unique) solution of (12)
 corresponding to the exact data $f_0, g_0 \in L^2(\R^2)$ and that $f,
 g \in L^2 (\R^2)$ are measured data satisfying $\| f - f_0\|_2 \leq
\varepsilon$ and $\| g - g_0\|_2 \leq \varepsilon$ where $\|.\|_2$ is the
 $L^2(\R^2)$-norm.

Then we can construct from $g, f$ a function $v_\varepsilon \in L^2$ such that
\begin{eqnarray}
\left\| {v_\varepsilon   - v_0 } \right\|_2  \le \sqrt {C\varepsilon ^{2 -
\gamma }  + \eta (\varepsilon )}
\nonumber
\end{eqnarray}
where $C$ is constant and $\eta(\varepsilon) \to 0$ as $\varepsilon
\downarrow 0$.

Moreover, if we assume in addition that $v_0  \in H^m (\R^2 )\cap L^1(
 \R^2), m > 0$ and \\ $0<\varepsilon<min\{e^{-e^2},e^{-4m^2}\}$ then
\begin{eqnarray}
\left\| {v_\varepsilon   - v_0 } \right\|_2  < D\left( {\ln \left( {\frac{1}
{\varepsilon }} \right)} \right)^{ - m}
\nonumber
\end{eqnarray}
where $D > 0$ depends on $v_0$.
}

{\bf Theorem 2}
{\sl

With $v_\varepsilon$ as in theorem 1, we have
\begin{eqnarray}
v_\varepsilon (x,t) = \sum\limits_{n = - \infty}^{+\infty} \sum\limits_{|m|
\leq |n|}v_\varepsilon\left(\frac{m\pi}{a_\varepsilon}, \frac{n\pi}
{a_\varepsilon}\right)S(m, \pi/a_\varepsilon)(x)S(n, \pi/a_\varepsilon)(t)
\nonumber
\end{eqnarray}
where
\begin{eqnarray}
S(p,d)(z) = \frac{sin[\pi(z - pd)/d]}{\pi(z - pd)/d},\,\,\,p \in \Z,
d > 0.
\nonumber
\end{eqnarray}
}

{\bf 3. Proofs}

{\bf Proof of theorem 1}

We put
\begin{eqnarray}
F(x,t) =   2R * f(x,t) - S * g(x,t) + 4f(x,t)
\nonumber
\end{eqnarray}
and
\begin{eqnarray}
F_0 (x,t) =   2R * f_0 (x,t) - S * g_0 (x,t) +4f_0(x,t)
\nonumber
\end{eqnarray}
then
\begin{eqnarray}
\left\| {\hat F - \hat F_0 } \right\|_2  &=& \left\| {F - F_0 } \right\|_2
\nonumber\\
&\le& (4 + 2\left\| R \right\|_1 )\left\| {f - f_0 } \right\|_2  + \left\| S
 \right\|_1 \left\| {g - g_0 } \right\|_2
\nonumber\\
&\le& \left( {4 + 2\left\| R \right\|_1  + \left\| S \right\|_1 } \right)
\varepsilon.
\nonumber
\end{eqnarray}

Put
\begin{eqnarray}
v_\varepsilon  (x,t) = \frac{1}{{2\pi }}\int\limits_{D_\varepsilon}
{\frac{{\hat F(z,r)}}{{\hat S (z,r)}}e^{i(xz + tr)}dzdr}
\end{eqnarray}
where $D_\varepsilon   = \left\{ {(z,r)/\left| z \right| \le b_\varepsilon  }
 \right.$ and $\left. {\left| r \right| \le b_\varepsilon ^2 } \right\}$ with
 $b_\varepsilon   = \frac{1}{{\sqrt 2 \sqrt {\sqrt 2  + 1} }}\ln \frac{4}
{{\varepsilon ^\gamma  }}$.

We have
\begin{eqnarray}
\left\| {v_\varepsilon   - v_0 } \right\|_2^2  = \left\| {\hat v_\varepsilon
 - \hat v_0 } \right\|_2^2  &=&\int\limits_{D_\varepsilon} {\left|
 {\frac{{\hat F(z,r) - \hat F_0 (z,r)}}{{\hat S(z,r)}}} \right|^2 dzdr}  +
\int\limits_{\R^2 \backslash D_\varepsilon} {\left| {\hat v_0 (z,r)}
 \right|^2 dzdr}
\nonumber\\
&\le& {\varepsilon ^{2 - \gamma } \left( {4 + 2\left\| R \right\|_1  +
\left\| S \right\|_1 } \right)^2 } + \int\limits_{\R^2 \backslash
D_\varepsilon} {\left| {\hat v_0 (z,r)} \right|^2 dzdr}
\nonumber
\end{eqnarray}

If we put $\eta (\varepsilon ) = \int\limits_{\R^2 \backslash
D_\varepsilon} {\left| {\hat v_0 (z,r)} \right|^2 dzdr}$ then $\eta
(\varepsilon ) \to 0$ as $\varepsilon  \downarrow 0$.

Now, we assume $v_0  \in H^m (\R^2 ), m > 0$, put
\begin{eqnarray}
a_\varepsilon   = \frac{{\sqrt 2 }}{{\sqrt {\sqrt 2  + 1} }}\ln \left(
{\frac{{1/\varepsilon }}{{\ln ^m (1/\varepsilon )}}} \right) > 1
\nonumber\\
Q_\varepsilon   = [ - a_\varepsilon  ,a_\varepsilon  ] \times [ -
a_\varepsilon  ,a_\varepsilon  ]
\nonumber
\end{eqnarray}
and
\begin{eqnarray}
v_\varepsilon  (x,t) = \frac{1}{{2\pi }}\int\limits_{Q_\varepsilon  }
{\frac{{\hat F(z,r)}}{{\hat S(z,r)}}e^{i(xz + tr)} dzdr}
\nonumber
\end{eqnarray}

We have
\begin{eqnarray}
\left\| {v_\varepsilon   - v_0 } \right\|_2^2  &=& \int\limits_{Q_\varepsilon
} {\frac{{\left| {\hat F(z,r) - \hat F_0 (z,r)} \right|^2 }}{{\left| {\hat
S(z,r)} \right|^2 }}dzdr}  + \int\limits_{\R^2 \backslash
Q_\varepsilon} {\left| {\hat v_0 (z,r)} \right|^2 dzdr}
\nonumber\\
 &\le& 4\varepsilon^2 (4 + 2\left\| R \right\|_1  + \left\| S \right\|_1 )^2
 e^{\sqrt 2 \sqrt {\sqrt 2  + 1} a_\varepsilon  }  + \int\limits_{\R^2
 \backslash Q_\varepsilon  } {\frac{{(z^2  + r^2 )^m \left| {\hat v_0 (z,r)}
 \right|^2 }}{{(z^2  + r^2 )^m }}dzdr}
\nonumber\\
 &\le& C_1 \left( {\varepsilon ^2 e^{\sqrt 2 \sqrt {\sqrt 2  + 1}
a_\varepsilon  }  + \frac{1}{{(2a_\varepsilon^2)^m }}} \right)
\nonumber
\end{eqnarray}
where
\begin{eqnarray}
C_1  = \max \left\{ {4\left( {4 + 2\left\| R \right\|_1  + \left\| S
\right\|_1 } \right)^2 ,\left\| {(z^2  + r^2 )^{m/2} \hat v_0 (z,r)}
\right\|_2^2 } 2^m\right\}
\nonumber
\end{eqnarray}

This implies that
\begin{eqnarray}
\left\| {v_\varepsilon   - v_0 } \right\|_2^2  &\le& C_1 \left[
{\varepsilon ^2 \left( {\frac{{1/\varepsilon }}{{\ln ^m (1/\varepsilon )}}}
\right)^2  + \frac{1}{{\ln ^{2m} \left( {\frac{{1/\varepsilon }}{{\ln ^m
 (1/\varepsilon )}}} \right)}}} \right]
\nonumber\\
&\le& C_1 \left[ {\frac{1}{{\ln ^{2m} (1/\varepsilon )}} +
\frac{1}{{\ln ^{2m} \left( {\frac{{1/\varepsilon }}{{\ln ^m
(1/\varepsilon)}}} \right)}}} \right]
\nonumber\\
&\le& C_1 \left[ {\frac{1}{{\ln ^{2m} (1/\varepsilon )}} + \frac{{2^m }}
{{\ln ^{2m} (1/\varepsilon )}}} \right] = D^2 \frac{1}{{\ln ^{2m} (1/
\varepsilon )}}
\nonumber
\end{eqnarray}
where $D = \sqrt {C_1 (1 + 2^m )} $.

This completes the proof.

{\bf Proof of theorem 2}

We have
\begin{eqnarray}
supp\,\,\, \hat{v_\varepsilon}\subset D_\varepsilon \subset
[-a_\varepsilon,a_\varepsilon] \times [-a_\varepsilon,a_\varepsilon].
\nonumber
\end{eqnarray}

As in [AGLT], p.121, we have theorem 2.

This completes the proof.

\bf 4. Numerical results\ \par
\ \ \ \rm We present some results of numerical comparison of the regularized
representation of the solution given by theorem 2 and the corresponding exact solution of the problem.\ \par
Let the problem\ \par
\ \par
\begin{equation}
\Delta u-{{\ds \partial  u}\over{\ds \partial  t}}
=0,\; \; (x,y)\in  \R  \times  (0,2),\; t>0
\end{equation}
\par
\begin{equation}
u(x,1,t)={{\ds 1}\over{\ds t}}  {\rm e}^{  {{
-x^{2}-1}\over{  4t}}  }  ;\; \; u(x,2,t)={{\ds 1}\over{\ds
t}}  {\rm e}^{  {{  -x^{2}-4}\over{  4t}} }  ;\; u(x,y,0)=0
\end{equation}
whose the unknown is \ \par
\par
\begin{equation}
v(x,t)=u(x,0,t),\quad x\geq x_0>0, \quad t>0
\end{equation}
The exact solution of this problem is \ \par
\par
$$ v(x,t)={{\ds 1}\over{\ds t }}  {\rm e}^{{{-x^{2}}\over{4t}} }  .$$
The approximated solution is calculated from the expansion of two-dimensional
Sinc series given by theorem 2 associated to formula (13) in which\ \par
\ \ \quad $ \hat F= \ds 4\frac{ {\rm e}^{-\sqrt{r^2+z^4}}} {\sqrt{r^2+z^4}}, \qquad
\hat S = 2{\rm e}^{-\sqrt{r^2+z^4}} $
\ \par

Thus we have \ \par
$$ {\left\Vert \hat F-\hat F_{0}  \right\Vert }_{ L^{2}(\R
^{2})}  =\varepsilon  $$
which is a perturbation similar to the one operated on the two given
functions $ f$ and $ g$.\ \par
With $ \varepsilon  ={{\ds 1}\over{\ds 50}}  ,\  N=50\
$(the size of the double series) and for $ (x,t)\in  [0.25,1.3]\times  [0,4]$
we have drawn the corresponding approximate surface solution $ (x,t)
\longrightarrow
v_{\varepsilon  }(x,t)$ in Fig.1.\ \par
To calculate the double integral in (10) we have used the rectangle
rule which gives good accuracy if one integrates on the interval $ [\varepsilon
,1/\varepsilon  ]\times  [\varepsilon  ,1/\varepsilon  ]$. The time of calculus
with a good computer is very long: 2 hours for 900 points $ M=(x,t)\in
[0.25,1.3]\times  [0,4]$. It is the reason for which we are limited ourselves
to a relatively small size of the double series ($ N=20)$. For comparison
in Fig.2 we have drawn the exact solution $ (x,t)\longrightarrow  v(x,t)$.\ \par

\centerline{
\includegraphics[width=3.5in]{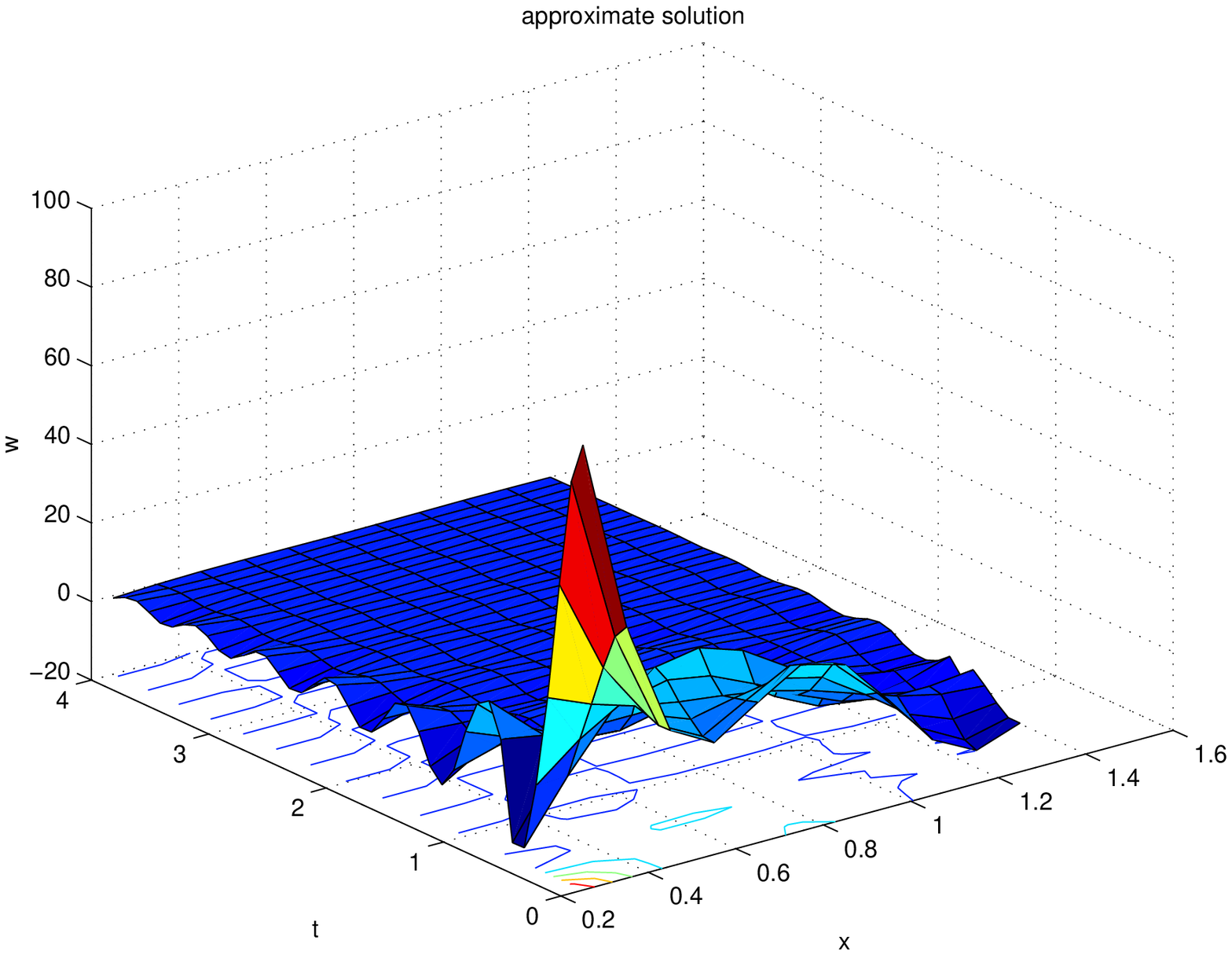}}
\centerline{Fig.1: regularized solution  of the problem (14),
(15)}\par \smallskip
\par
\centerline{
\includegraphics[width=3.5in]{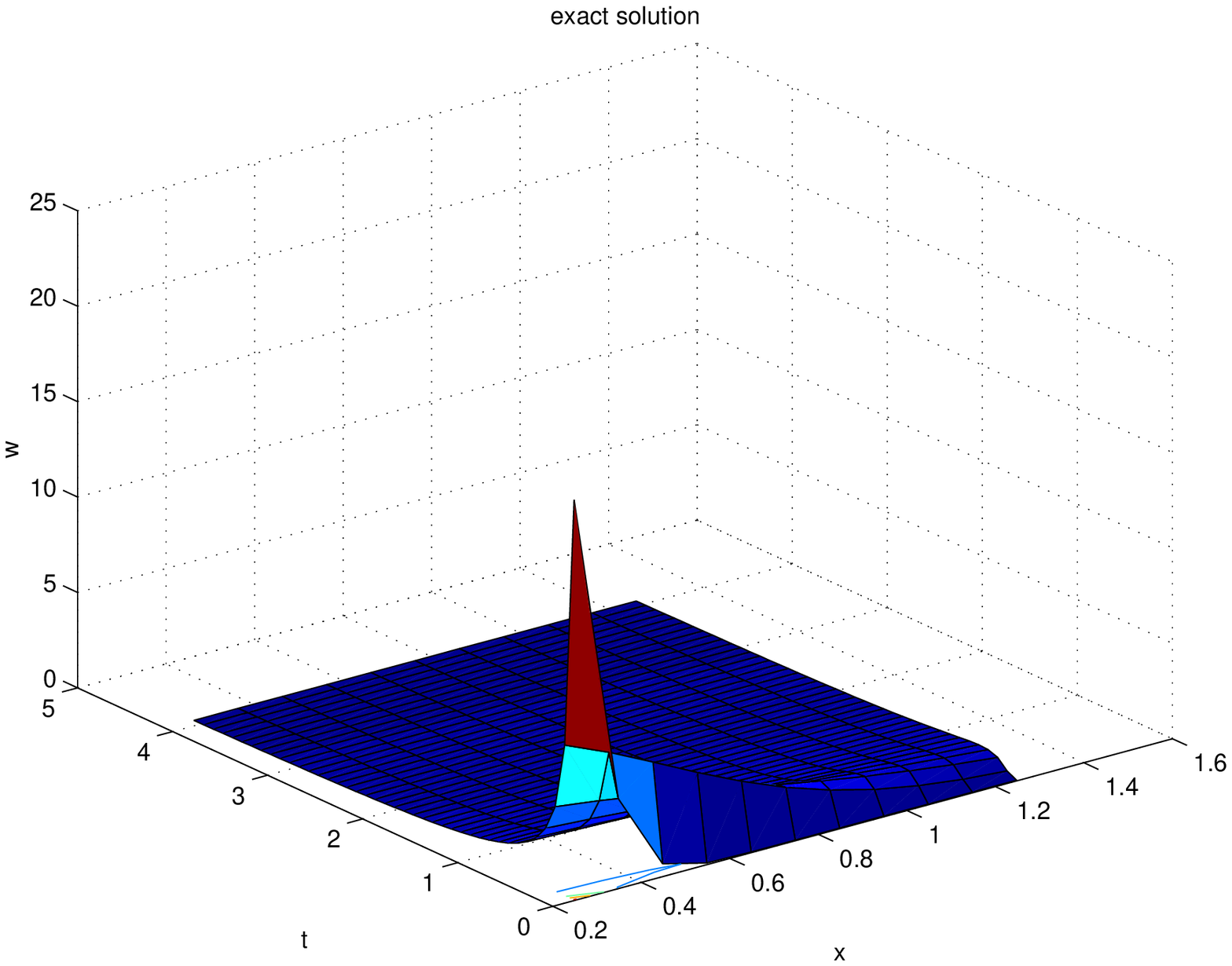}}
\centerline{Fig.1: exact solution  of the problem (14), (15)}\par
\smallskip
\par

Using the same method as previously we have drawn in Fig.3 the
surface $ (x,t)\in  [0,1]\times  [0,4]\longrightarrow v_{\varepsilon
}(x,t)$ which is the regularization of the following problem\ \par

\begin{equation}
 \Delta u-{{\ds \partial  u}\over{\ds \partial  t}}
=0,\; \; (x,y)\in  \R  \times  (0,2),\; t>0
\label{}
\end{equation}

\begin{equation}
u(x,1,t)=0;\; \; u(x,2,t)={{\ds 1}\over{\ds t}}  {\rm e}^{
{{  -x^{2}-4}\over{  4t}}  } ;\; u(x,y,0)=0
\label{}
\end{equation}
the unknown being $ v(x,t)=u(x,0,t)$. The exact solution \ \par
$$ v(x,t)=-{{\ds 1}\over{\ds t }}  {\rm e}^{
{{  -x^{2}-4}\over{4t}}  }  .$$
is represented in Fig.4.\ \par

\ \ \ \ \par
\centerline{
\includegraphics[width=3.5in]{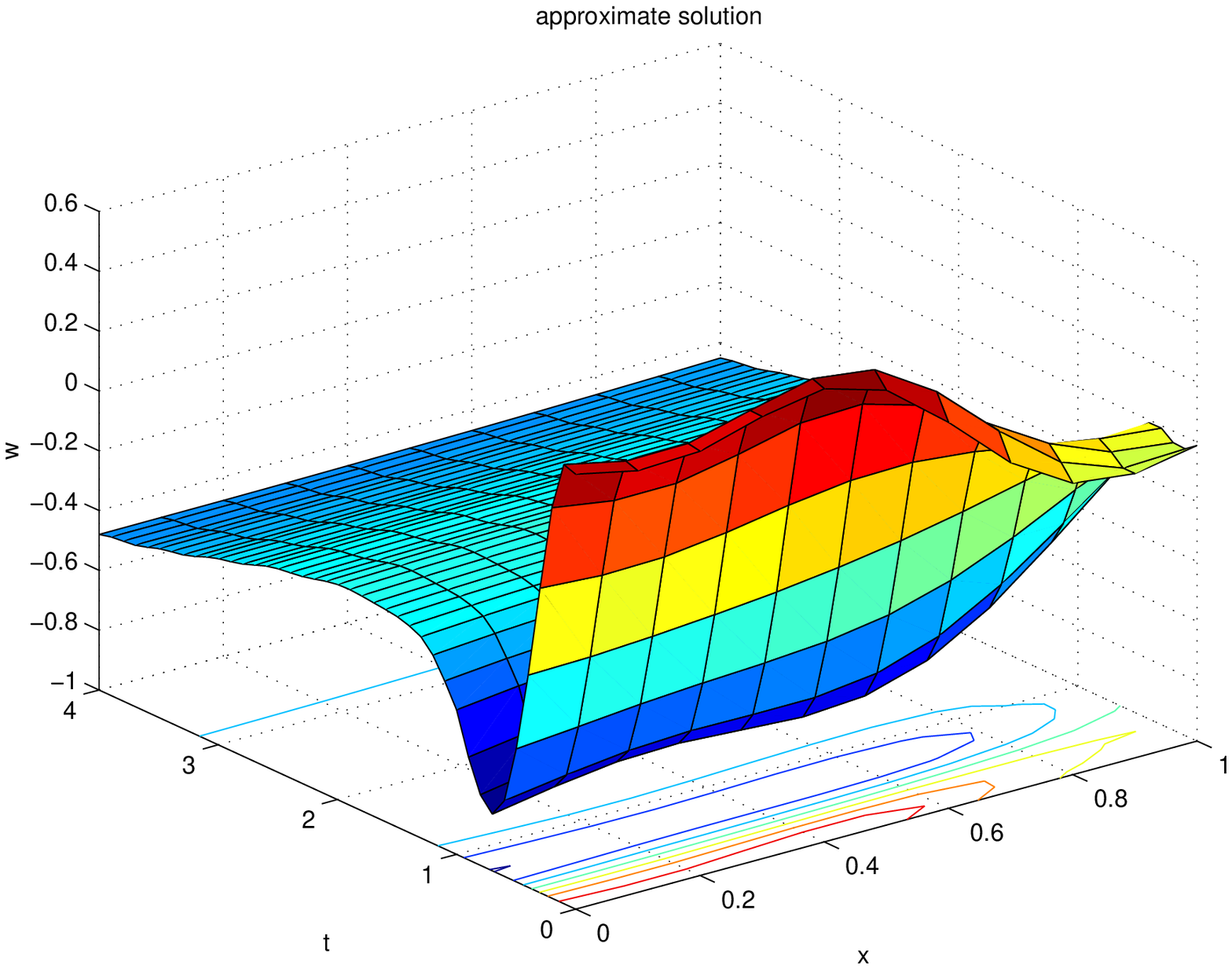}}
\centerline{Fig.3: regularized solution  of the problem (17),
(18)}\par
\smallskip
\par
\centerline{
\includegraphics[width=3.5in]{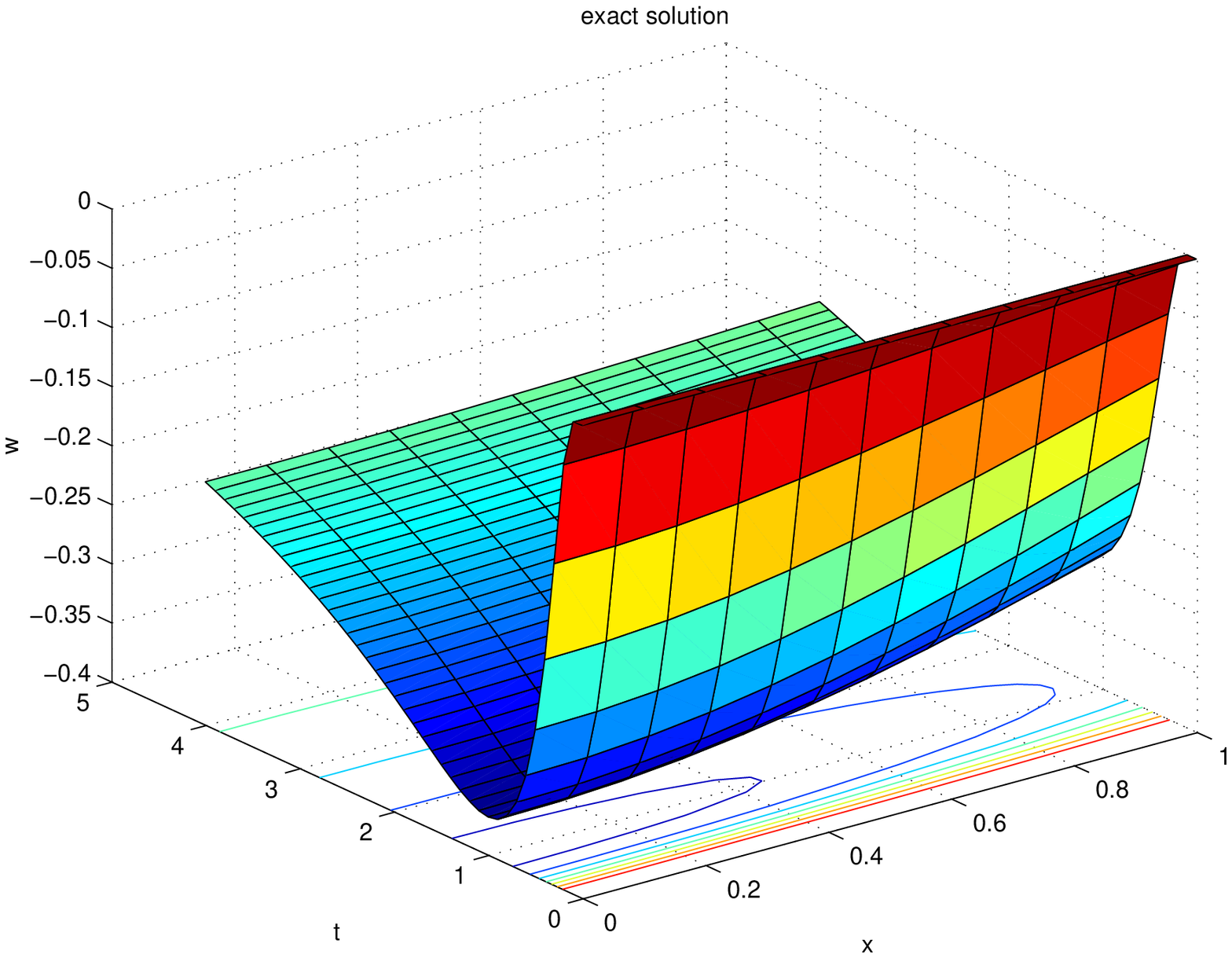}}
\centerline{Fig.4: exact solution  of the problem (17), (18)}\par
\smallskip
\par

\begin{center}
REFERENCES
\end{center}

[AGVT] D. D. Ang, R. Gorenflo, L. K. Vy and D. D. Trong.
{\sl Moment theory and some Inverse Problems in Potential Theory and Heat
 Conduction}, Lecture Notes in Mathematics 1792, Springer, (2002).

[B] J. V. Beck, B. Blackwell and C. R. St. Clair, Jr., {\sl Inverse Heat
Conduction, Ill-posed Problem}, Wiley, New York, (1985).

[C] A. Carasso, {\sl Determining surface temperatures from interior
observations}, SIAM J. Appl. Math. {\bf 42} (1981), 558-547.

[E] Erdelyi et al., {\sl Tables of Integral Transforms}, Vol. 1, Mc
Graw-Hill, New York, (1954).

[EM] H. Engl and P. Manselli, {\sl Stability estimates and regularization for
 an inverse heat conduction problem}, Numer. Funct. Anal. and Optim. {\bf 10}
 (1989), 517-540.

[F] A. Friedman, {\sl Partial Differential Equations of Parabolic Type},
 Englewood Cliff., N. J., (1964).

[LN] T. T. Le and M. P. Navarro, {\sl Surface Temperature From Borehole
Measurements : Regularization and Error Estimates}, Internl. J. Math and Math
 Sci. (1995).

[S] Stenger Fr., {\sl Numerical methods based on Sinc and analytic functions}
, Springer Verlag, Berlin - Heidelberg - New York, (1993).

[TA] A. N. Tikhonov and V. Y. Arsenin, {\sl Solutions of ill-posed problems},
 Winston, Washington, (1977).

[TV] G. Talenti and S. Vessella, {\sl Note on an ill-posed problem for the
heat equation}, J Austral. Math. Soc. {\bf 32} (1981), 358-368.

\end{document}